\documentclass{article}
\usepackage{amsmath,amsfonts,amssymb,theorem,verbatim}
\usepackage{amsfonts,amssymb,graphics,theorem}

\newtheorem{lem}{Lemma} 
\newtheorem{Thm}{Theorem}
 \newtheorem{Lemma}{Lemma}
\newtheorem{Prop}{Proposition} 
{\theorembodyfont{\normalfont}

}

\def\BEN{\begin{enumerate}}  \def\BI{\begin{itemize}}
\def\EEN{\end{enumerate}}   \def\EI{\end{itemize}} 
   
  \def\no{\noindent}

\def\beq{\begin{eqnarray}} \def\eeq{\end{eqnarray}}

\def\al*#1{\begin{align*}#1\end{align*}}

\def\ga*#1{\begin{gather*}#1\end{gather*}}

\def\alat*#1#2{\begin{alignat*}{#1}#2\end{alignat*}}
\def\bea{\begin{eqnarray*}}
\def\eea{\end{eqnarray*}}
\def\ml*#1{\begin{multline*}#1\end{multline*}}

  \def\ovl{\overline}

\def\mc{\mathcal}

\def\le{\left} \def\ri{\right} \def\i{\infty}

  \def\td{\text{\rm d}}

 \def\WT{\widetilde}

\def\a{\alpha} \def\b{\beta}
  \def\d{\delta}   

\def\e{\epsilon} 
\def\m{\mu} 
     \def\s{\sigma}

\newcommand{\proof}{\no{\it Proof\ }}
\newcommand{\halmos}{\vspace{3mm} \hfill $\Box$}
\newcommand{\exit}{{\mbox{\, \vspace{3mm}}} \hfill\mbox{$\square$}}

\begin{document}
\title{The probability of exceeding a piecewise deterministic barrier
by the heavy-tailed renewal compound process}

\author{Zbigniew
Palmowski\footnote{University of Wroc\l aw, pl. Grunwaldzki 2/4,
50-384 Wroc\l aw, Poland, E-mail: zpalma@math.uni.wroc.pl}
 \hspace{0.4cm}  Martijn Pistorius\footnote{Department
of Mathematics, King's College London, Strand, London WC2R 2LS,
UK, E-mail:Martijn.Pistorius@kcl.ac.uk}
 } \maketitle

 \begin{abstract}{\small\no
We analyze the asymptotics of crossing a high piecewise linear
barriers by a renewal compound process with the subexponential
jumps. The study is motivated by ruin probabilities of two
insurance companies (or two branches of the same company) that
divide between them both claims  and premia in some specified
proportions when the initial reserves of both companies tend to
infinity.}
\end{abstract}
\medskip


\thispagestyle{empty}
\addtocounter{page}{-1}
\newpage

\section{Motivation}

The study of boundary crossing probabilities of a stochastic
process with heavy tailed increments has applications in fields
such as queuing theory, insurance and finance. In this paper we
consider the problem of a renewal process crossing a high
piecewise linear boundary. This study is in particular motivated
by ruin probabilities of two insurance companies with proportional
claims (see \cite{APP}) and the steady state distribution of a
tandem queue with two servers (see \cite{Mandjes}). To be more
precise, we set
\begin{equation}\label{St} S_t=\sum_{i=1}^{N_t}
\sigma_i\end{equation} for $N_t$ a renewal process with i.i.d.
inter-arrival times $\zeta_i$, and the claims $\sigma_i$ are
i.i.d. random variables independent of $N(t)$, with the
distribution function $F(x)$. We shall denote by $\lambda$ and
$\m$ the reciprocals of the means of $\zeta_i$  and $\sigma_i$,
respectively. Let the boundaries $b_1, b_2$ given by
$$
b_1(t) = b_1(t;x_1) = x_1 + p_1t,\quad\quad b_2(t) = b_2(t;x_2) =
x_2 + p_2t,
$$
where we assume that
\begin{equation}\label{p1p2}
p_1 > p_2,\quad\quad\quad
p_2>\rho:=\frac{\lambda}{\m}=E[\sigma]/E[\zeta]
\end{equation}
and consider the following boundary crossing probabilities:
\begin{eqnarray*}
\psi_{\wedge}(x_1,x_2) &=& P(\exists t\ge 0: S_t > b_1(t)\wedge
b_2(t))\\ \psi_{\vee}(x_1,x_2) &=& P(\exists t\ge 0: S_t > b_1(t)\vee b_2(t))\\
\psi_{\times}(x_1,x_2) &=& P(\exists t\ge0, u\ge0: S_t > b_1(t),
S_u
> b_2(u)),
\end{eqnarray*}
where $x\vee y = \max\{x,y\}$ and $x\wedge y=\min\{x,y\}$.

The $\psi_{\wedge}(x_1,x_2)$ describes the ruin probability of at
least one insurance company. $\psi_{\vee}(x_1,x_2)$ corresponds to
the simultaneous ruin of the insurance companies. Finally,
$\psi_{\times}(x_1,x_2)$ describes probability that both companies
will have ruin. First assumption in (\ref{p1p2}) means that the
second company receives less premium per amount paid out and
second one is the stability condition under which reserves of both
insurance companies tend to infinity. The solutions of the
"degenerate two-di\-men\-si\-o\-nal" ruin problems strongly depend
on the relative position of the vector of premium rates $p=(p_1,
p_2)$ with respect to the proportions vector $(1,1)$. Namely, if
the initial reserves satisfy $x_2 \leq x_1$, the two lines do not
intersect. It follows therefore that the barriers the  "$\wedge$"
and "$\vee$" ruin always happen for the second and first company
respectively. In this case the asymptotics follows from
one-dimensional ruin theory  -- see e.g. Rolski et al. (1999).
Therefore we focus here on the opposite case when $x_1<x_2$.

In this paper we derive the exact first order asymptotics of these
ruin probabilities if $x_1, x_2$ tend to infinity according to a
ray in the positive quadrant (i.e. $x_1/x_2$ is constant) if the
claims follow a subexponential distribution. We model the claims
by subexponential distributions since many catastrophic events
like earthquakes, storms, terrorist attacks etc are used in their
description. Insurance companies use e.g. the lognormal
distribution (which is subexponential) to model car claims -- see
Rolski et al. (1999) or Embrechts et al. (1997) for the further
background.

\section{Main result}
In order to state our results we start with recalling some
notions. A distribution function $G$ on $[0,\i)$ is subexponential
($G\in \mathcal{S}$) if
$$\lim_{x\to\infty}\overline{G^{*2}}(x)/\overline{G}(x)=2\;,$$ where
$G^{*2}$ is the fold of $G$ with itself. The integrated tail
distribution $G_I$ of $G$ is defined as $G_I(x) = \nu^{-1}\int_0^x
\ovl G(y)\td y$, where $\nu = \int_0^\infty \overline{G}(x) \,{\rm
d}x$. In fact we will impose the stronger condition, that is,
\begin{equation}\label{RV}
\overline{G}(x)\sim Ax^{-\alpha},\quad 1<\alpha<2.
\end{equation}
Note that $G$ is long-tailed: \begin{equation}
  \label{eq:1}
  \text{$\overline{G}(x)>0$ for all $x$},\qquad
  \lim_{x\to\infty}\frac{\overline{G}(x-h)}{\overline{G}(x)} = 1,
  \quad\text{for all fixed $h>0$}.
\end{equation}

 Let $F$ denote
the distribution of the claim size $\sigma$ and define
\begin{equation}\label{H}H(aK,K)= \int_0^\infty
\overline{F}\left(\max\left\{aK+m_1 t,K+m_2 t\right\}\right) \;\td
t
\end{equation}
and \begin{equation}\label{J}J(aK,K)= \int_0^\infty
\overline{F}\left(\min\left\{aK+m_1 t,K+m_2 t\right\}\right) \;\td
t
\end{equation}
 for $m_i = p_iE[\zeta] - E[\sigma]$ ($i=1,2$).

In view of (\ref{p1p2}) we note that if $a \ge 1$, $b_1\vee b_2 =
b_1$ and $b_1\wedge b_2=b_2$, so that the crossing problems are
one-dimensional. If $F_I$ follows a subexponential distribution,
Veraverbeke's theorem (see Veraverbeke (1977), Embrechts and
Veraverbeke (1982) and also Zachary (2004) for a short proof)
implies that
$$\psi_{\wedge}(aK,K) \sim
\frac{1}{m_2 \m}\overline{F}_I(K), \quad \psi_{\vee}(aK,K)\sim
\frac{1}{m_1 \m}{ \overline{F}_I(aK)},\quad K\to\i,
$$ where we write $f(x)\sim g(x)$ $(x\to\i)$ if $\lim_{x\to\i}f(x)/g(x) =
1$.
In the opposite case that $a<1$ the asymptotics of
$\psi_{\wedge}/\psi_{\vee}/\psi_{\times}$ are as follows:
\begin{Thm}\label{thm}
Let $a<1$. Assume that $E[\zeta^{3}]<\i$. If $F$ satisfies
(\ref{RV}) with $1<\alpha<2$, then it holds that, as $K\to\infty$,
\begin{eqnarray}\label{eq:asymppr}
\psi_{\wedge}(aK,K)&\sim&\frac{1}{m_1\m} \ovl F_I(aK) +
\frac{1}{m_2\m}\overline{F}_I(K) - H(aK,K)\sim J(aK,K)\;,\\
\psi_{\vee}(aK,K)&\sim& \psi_{\times}(aK,K) \sim H(aK,K).
\end{eqnarray}
\end{Thm}


It is worth observing that, in contrast to the case of light-tails
(see \cite[Thm. 1]{APP}), the asymptotic probability of crossing
both $b_1$ and $b_2$ at the same time appears in the asymptotics
for $\psi_{\wedge}(aK,K)$ (and is equal to $H(aK, K)$). In the
setting of (re)insurance companies with proportional claims (see
\cite{APP} for details) this result agrees with what we might
expect - "large" claims cause often bankruptcy not only of the
insurance company but of a whole chain of reinsurers.

The proof of the main result follows from the Lemmas ... The paper
is organized as follows. In Section \ref{keylemmas} we present
useful facts and in Section \ref{mainestimates} main (lower and
upper) bounds giving the Theorem \ref{thm}.

\section{Key Lemmas}\label{keylemmas}

Let $T_0=0, T_n=\sum_{i=1}^n\zeta_i$ and $\Xi_0=0,
\Xi_n=\sum_{i=1}^n\sigma_i$ denote the two random walks of claims
and interarrival times and
define the associated random walks $(S^{(1)}_n)_{n\ge 0}$
and $(S^{(2)}_n)_{n\ge 0}$ by
$$
S^{(i)}_n = \Xi_n - p_i T_n\quad i=1,2.
$$
Let $M^{(i)}_n = \max_{m\leq n} S^{(i)}_m$ and
$M^{(i)}_{[n,\i)} = \max_{m\ge n} S^{(i)}_m$
with $M^{(i)}_\i = M^{(i)}_{[0,\i)}$.
Also let
$$
T = \frac{x_2-x_1}{m_1-m_2} = \frac{x_2-x_1}{(p_1-p_2)E[\zeta]}
$$
be the epoch that the lines $\WT b_1(t) = x_1 + m_1 t$ and $\WT
b_2(t) = x_2 + m_2 t$ cross and let
$$
\WT T =  \frac{x_2-x_1}{p_1-p_2} = \frac{T}{\lambda}
$$
be the epoch that the lines $b_1(t) = x_1 + p_1 t$ and $
b_2(t) = x_2 + p_2 t$ cross.
We denote
$$c=\frac{1-a}{m_1-m_2}$$
and note that $T=cK$ if $x_1=aK$ and $x_2=K$.
Define
$$
N^* = \min\{n\ge 1: T_n > \WT T\} - 1
$$
and note that $N^*$ is a stopping time w.r.t. $\s(\{\Xi_n,
T_n\}_n)$. Further, in view of the strong law of large numbers
$N^*/T \to 1$. From the Theorem V.5.14 of Petrov (1972) we have
the following upper bound:
\begin{lem}\label{uuppbound}
Let $\zeta^3<\infty$. Denote $\rho=E|\zeta-E\zeta|^3/s^2$ for
$s^2=E(\zeta-E\zeta)^2$. Then
$$\left|P\left(\frac{T_n-nE\zeta}{\sqrt{n}}<x\right)-\Phi(x)\right|\leq
A\frac{\rho}{\sqrt{n}(1+|x|)^3}$$ for all $x$ and fixed $A$.
\end{lem}
Now, we have
\begin{eqnarray*}\lefteqn{P(N^*>(1+\e)T)=P(T_{[(1+\e)T]}<\tilde{T})}\\&&=
P\left(\frac{T_{[(1+\e)T]}-(1+\e)T\frac{1}{\lambda}}{\sqrt{T(1+\e)}}<-\sqrt{T}\frac{\e}{\lambda\sqrt{(1+\e)}}\right)\leq
\Phi\left(-\sqrt{T}\frac{\e}{\lambda\sqrt{(1+\e)}}\right)\\&&\ \ \
\ +A\frac{\rho(1+\e)\lambda^3}{T^2\e^3}= {\rm
O}(T^{-2})=o(T^{-\alpha}).
\end{eqnarray*}
Hence
\begin{equation}\label{CLT}
P(N^*>(1+\e)T)={\rm o}(\overline{F}(T)),
\end{equation}
Similarly,
\begin{equation}\label{CLT2}
P(N^*<(1-\e)T)={\rm o}(\overline{F}(T)).
\end{equation}
and
\begin{eqnarray}\label{CLT3}
P(T_T<(1-\e)E\zeta T)&=&{\rm o}(\overline{F}(T)),\\
P(T_T>(1+\e)E\zeta T)&=&{\rm o}(\overline{F}(T)).
\end{eqnarray}

 A key step is to note that the ruin probabilities can be
defined in terms of these quantities as follows:

\begin{Prop} It holds that
\begin{eqnarray*}
\psi_{\vee}(aK,K) &=& P(M^{(2)}_{N^*} > K) + P(M^{(2)}_{N^*} \leq
K,
M^{(1)}_{[N^* + 1,\i)} > aK),\\
\psi_{\wedge}(aK,K) &=& P(M^{(1)}_{N^*} > aK) + P(M^{(1)}_{N^*}
\leq aK,
M^{(2)}_{[N^* + 1,\i)} > K),\\
\psi_{\times}(aK,K) &=& P(M_\i^{(2)} > K) - P(M_\i^{(2)} > K, M_\i^{(1)} \leq aK)\\
&=& P(M_\i^{(2)} > K) + P(M^{(1)}_{N^*} \leq aK,
M^{(1)}_{[N^*+1,\i)} > aK)\\&& - P(M^{(1)}_{N^*} \leq aK,
M^{(2)}_{[N^*+1,\i)} > K).
\end{eqnarray*}
\end{Prop}

Indeed note that
\begin{eqnarray*}
\lefteqn{P(M_\i^{(2)} > K, M_\i^{(1)} \leq aK) =
P(\max\{M^{(2)}_{N^*},M^{(2)}_{[N^*+1,\i)} \}>K, M^{(1)}_{N^*}
\leq aK,
M^{(1)}_{[N^*+1,\i)} \leq aK)}\\
&&= P(M^{(2)}_{[N^*+1,\i)} > K, M^{(1)}_{N^*} \leq aK,
M^{(1)}_{[N^*+1,\i)} \leq aK)\\
&&= P(M^{(1)}_{N^*} \leq aK, M^{(2)}_{[N^*+1,\i)} > K) -
P(M^{(1)}_{N^*} \leq aK, M^{(1)}_{[N^*+1,\i)} > aK).
\end{eqnarray*}

Note that in view of Veraverbeke's theorem
$$
P(M_\i^{(2)} > K) \sim \frac{1}{m_2\mu}\ovl F_I(K).
$$
To prove the Theorem \ref{thm} we will estimate the other terms.

\section{Main estimates}\label{mainestimates}
Before we proceed we first introduce some extra notation and
collect some auxiliary results.
Let
\begin{eqnarray}\label{innereprH}
H(x,y) &=& I_{[0,T]}^{(2)}(x)+I_{[T,\infty]}^{(1)}(y), \\
J(x,y) &=& I_{[0,T]}^{(1)}(x)+I_{[T,\infty]}^{(2)}(y),
\label{innereprJ}
\end{eqnarray}
where $I_{[x,y]}^{(i)}(w)=\int_x^y \overline{F}\left(w+m_i
t\right) \;\td t$ for $i=1,2$. Note that definition of $H(ak,k)$
and $J(aK,K)$ agree with (\ref{H}) and (\ref{J}) respectively.

The following result relates
$H(aK,K)$ to $I^{(2)}_{[0,\i]}(aK)$ and $J(aK,K)$ to
$I^{(1)}_{[0,\i]}(K)$:
\begin{Lemma}\label{lem:HI2}
The following hold true: \BI
\item[(i)] $I^{(1)}_{[T,\i]}(aK) = (m_2/m_1)I^{(2)}_{[T,\i]}(K) $,
\item[(ii)] $
H(aK,K) \leq
I^{(2)}_{[0,\infty]}(aK)$. %
\item[(iii)] $
J(aK,K) \leq (m_1/m_2) I^{(1)}_{[0,\infty]}(aK)$,
\item[(iv)] For $i=1,2$ it holds that $I^{(i)}_{[\a,\b]}(\zeta)
\leq \frac{\b-\a}{\a} I^{(i)}_{[0,\a]}(\zeta)$ where
$\a,\b,\zeta>0$ are constants with $\a<\b$. In particular,
$I^{(1)}_{[T,T(1+\b)]}(aK) \leq \b I^{(1)}_{[0,T]}(aK)$ and
$I^{(2)}_{[T,T(1+\b)]}(K) \leq \b I^{(2)}_{[0,T]}(K)$. \EI
\end{Lemma}
{\bf Proof:} (i) follows by a straightforward calculation. The
bounds in (ii) and (iii) follow then in view of (i), equations
(\ref{innereprH}) and (\ref{innereprJ}) and since $m_2\leq m_1$.

(iv) Since $\ovl F$ is decreasing it follows that
$I^{(i)}_{[\a,\b]}(\zeta)$ is bounded above by $(\b-\a)\ovl
F(\zeta + m_i\a)$ and $I^{(i)}_{[0,\a]}(\zeta)$ is bounded below
by $\a\ovl F(\zeta+m_i\a)$. Combining these two estimates proves
the statement.\exit

Since $\psi_\times \ge \psi_{\vee}$, we
note that it will be sufficient to prove the following estimates:

\begin{eqnarray}\label{prop1.1}
\liminf_{K\to\infty}\frac{\psi_{\wedge}(aK,K)}{J(aK,K)}&\ge& 1,\\
\liminf_{K\to\infty}\frac{\psi_{\vee}(aK,K)}{H(aK,K)}\ge 1,
\label{prop1.2}
\end{eqnarray}
and
\begin{eqnarray}
\limsup_{K\to\infty}\frac{\psi_{\wedge}(aK,K)}{J(aK,K)}&\leq&1,
\label{prop2.1}
\\
\limsup_{K\to\infty}\frac{\psi_{\times}(aK,K)}{H(aK,K)}\leq 1.
\label{prop2.2}
\end{eqnarray}

\begin{Lemma} Suppose that $F$ satisfies (\ref{RV}).
Given $\eta>0$, it holds for $K$ large enough that
\begin{eqnarray*}
P(M^{(j)}_{N^*} \leq a_jK, M^{(i)}_{[N^*+1,\i)} > a_iK) &\sim&
I^{(i)}_{[T,\i]}(a_iK).
\end{eqnarray*}
\end{Lemma}

\proof It follows from Lemmas \ref{lowerstad}, \ref{upperstad},
\ref{lowermiesz} and \ref{uppermiesz}. \halmos

\begin{Lemma} 
Suppose $F_I$ is long-tailed. Fix $0<\eta<1$. For $K$ large enough
it holds that
\begin{eqnarray*}
P(M^{(i)}_{N^*} > a_iK) &\sim& I^{(i)}_{[0,T]}(a_iK), \quad i=1,2.
\end{eqnarray*}
\end{Lemma}

\proof By Theorem 2 of Foss et al. (2006) for sufficiently large
$K$ and small $\e$ we have,
\begin{eqnarray*}
\lefteqn{P(M^{(i)}_{N^*} > a_iK) \geq P(M^{(i)}_{[T(1-\e)]} >
a_iK, N^*>T(1-\e))}\\&& \geq P(M^{(i)}_{[T(1-\e)]} > a_iK)-P(
N^*<T(1-\e)) = P(M^{(i)}_{[T(1-\e)]} > a_iK)-{\rm
o}(I^{(i)}_{[0,T]}(a_iK))\\&&\geq
(1-\eta/2)I^{(i)}_{[0,T(1-\e)]}(a_iK) =
(1-\eta/2)I^{(i)}_{[0,T]}(a_iK)-
(1-\eta/2)I^{(i)}_{[T(1-\e),T]}(a_iK)\\&&\geq
(1-\eta/2)I^{(i)}_{[0,T]}(a_iK)- (1-\eta/2)\e
I^{(i)}_{[0,T]}(a_iK)\geq (1-\eta)I^{(i)}_{[0,T]}(a_iK).
\end{eqnarray*}
Similarly,
\begin{eqnarray*}
\lefteqn{P(M^{(i)}_{N^*} > a_iK) \leq P(M^{(i)}_{[T(1-\e)]} >
a_iK, N^*<T(1+\e))+P(N^*>T(1+\e))}\\&& \leq P(M^{(i)}_{[T(1+\e)]}
> a_iK)+P( N^*>T(1+\e)) = P(M^{(i)}_{[T(1+\e)]} > a_iK)+{\rm
o}(I^{(i)}_{[0,T]}(a_iK))\\&&\leq
(1+\eta/2)I^{(i)}_{[0,T(1+\e)]}(a_iK) =
(1+\eta/2)I^{(i)}_{[0,T]}(a_iK)+
(1+\eta/2)I^{(i)}_{[T,T(1+\e)]}(a_iK)\\&&\leq
(1+\eta/2)I^{(i)}_{[0,T]}(a_iK)+ (1+\eta/2)\e
I^{(i)}_{[0,T]}(a_iK)\leq (1+\eta)I^{(i)}_{[0,T]}(a_iK).
\end{eqnarray*}

\halmos

\subsection{Lower bounds}
 Let
$$a_i=\left\{\begin{array}{lr}
a&\quad \mbox{for $i=1$,}\\
1&\quad \mbox{for $i=2$.}
\end{array}\right.$$

\begin{Lemma} Suppose $F_I\in\mc S$.
Given $\eta>0$, it holds for $K$ large enough that
\begin{eqnarray*}
P(M^{(i)}_{[N^*+1,\i)} > a_iK) &\ge&
(1-\eta)I^{(i)}_{[T,\i]}(a_iK) - \eta I^{(i)}_{[0,T]}(a_iK).
\end{eqnarray*}
\end{Lemma}

\proof Fix $1>\e>0$ and $0<\d<\e$. By the law of large numbers it holds for
$K$ large enough that
\begin{equation}\label{eq:lln}
P(S^{(i)}_{N^*} > -m_i( 1+ \e)N^*, N^* < T(1+\e)) > 1-\d.
\end{equation}
Denote by $\WT M^{(i)}$ an independent copy of $M^{(i)}$. It
follows by an application of Veraverbeke's theorem and Lemma
\ref{lem:HI2} that for $a_iK$ large enough it holds that
\begin{eqnarray*}
\lefteqn{P(M^{(i)}_{[N^*+1,\i)} > a_iK)}\\
&\ge& P(M^{(i)}_{[N^*+1,\i)} > a_iK, S^{(i)}_{N^*} > -(m_i + \e)N^*, N^* < T(1+\e))\\
&\ge& P(\WT M^{(i)}_{[0,\i)} > a_iK + m_i(1+\e)^2T)
P(S^{(i)}_{N^*} > -(m_i + \e)N^*, N^* < T(1+\e))\\
&\ge& (1-\d) I^{(i)}_{[0,\i]}(a_iK + m_i(1 + 3\e)T)\\
&\ge& (1-2\d)\int_0^\i \ovl F(a_iK + m_i(1+3\e)T + m_i t)\td t\\
&\ge& (1-2\d)\int_T^\i \ovl F(a_iK + m_i(1+3\e)t)\td t\\
&=& (1-2\d)\frac{m_i}{m_i(1 + 3\e)} I^{(i)}_{[T(1+ 3\e), \i)}(a_iK)\\
&\ge& (1-2\d)\frac{1}{1 + 3\e} \le(I^{(i)}_{[T,\i]}(a_iK) - 3\e
I^{(i)}_{[0,T]}(a_iK)\ri).
\end{eqnarray*}
\exit

\begin{Lemma} Suppose $F_I\in\mc S$.
For $i,j=1,2$, given $\eta>0$, it holds for $K$ large enough that
\begin{eqnarray*}
P(M^{(j)}_{N^*} \leq a_jK,
M^{(i)}_{[N^*+1,\i)} > a_iK) &\ge& I^{(i)}_{[T,\i]}(a_iK) - \eta I^{(i)}_{[0,T]}(a_iK).\\
\end{eqnarray*}
\end{Lemma}

\proof As the previous Lemma, replacing \eqref{eq:lln} by
$$
P(M^{(j)}_{N^*} \leq a_jK, S^{(i)}_{N^*} > -(m_i + \e)N^*, N^* <
T(1+\e)) > 1-\d.
$$
\exit

\begin{Lemma}\label{lowerstad} Suppose that $F$ satisfies (\ref{RV}).
Given $\eta>0$, it holds for $K$ large enough that
\begin{eqnarray*}
P(M^{(i)}_{[N^*+1,\i)} > a_iK) \geq (1-\eta) T
\overline{F}(a_iK+m_iT)+(1-\eta)I^{(i)}_{[T,\i]}(a_iK)
\end{eqnarray*}
\end{Lemma}

\proof Note that $$ P(M^{(i)}_{[N^*+1,\i)}
> a_iK)=P(S^{(i)}_{N^*+1} +\tilde{M}^{(i)}_{[0,\i)}> a_iK),$$
where $\tilde{M}^{(i)}_{[0,\i)}$ is all-time maximum from the
independent random walk distributed as $S^{(i)}_n$. Moreover, for
any $\epsilon>0$,
\begin{eqnarray}\lefteqn{P(S^{(i)}_{N^*+1} +\tilde{M}^{(i)}_{[0,\i)}> a_iK)\geq
P(\Xi_{N^*+1} +\tilde{M}^{(i)}_{[0,\i)}> a_iK +p_iT_{N^*+1},
T(1+\e)>N^*>T(1-\epsilon))}\nonumber\\&&\geq
P(\Xi_{[T(1-\epsilon)+1]} +\tilde{M}^{(i)}_{[0,\i)}> a_iK
+p_iT_{[T(1+\epsilon)+1]})\nonumber\\&& \geq
P(\Xi_{[T(1-\epsilon)+1]} +\tilde{M}^{(i)}_{[0,\i)}> a_iK
+p_iT(1+\epsilon)^2E\zeta)
-P(N^*<T(1-\epsilon))\nonumber\\&&-P(N^*>T(1+\epsilon))-
P(T_{[T(1+\epsilon)+1]}>(1+\e)^2E\zeta T )\nonumber\\&& \geq
P(\Xi_{[T(1-\epsilon)+1]} +\tilde{M}^{(i)}_{[0,\i)}> a_iK
+p_iT(1+\epsilon)^2E\zeta)- {\rm
o}(I_{[T,\infty]}^{(i)}(a_iK)).\label{normalref}
\end{eqnarray}
Further, for $1<\alpha< 2$ and
$$C_\alpha=\frac{1-\alpha}{\Gamma(2-\alpha)\cos(\pi\alpha/2)}$$
we derive
\begin{eqnarray*}
\lefteqn{P(\Xi_{[T(1-\epsilon)+1]} +\tilde{M}^{(i)}_{[0,\i)}> a_iK
+p_iT(1+\epsilon)^2E\zeta)}\\&&=
P\left(\frac{\Xi_{[T(1-\epsilon)+1]}-E\sigma
T(1-\epsilon)}{(A/C_\alpha)^{1/\alpha}T^{1/\alpha}(1-\epsilon)^{1/\alpha}}(A/C_\alpha)^{1/\alpha}T^{1/\alpha}(1-\epsilon)^{1/\alpha}
+\tilde{M}^{(i)}_{[0,\i)}> (a_i+\kappa)K +m_iT \right),
\end{eqnarray*}
where $\kappa=\epsilon (E\sigma-p_iE\zeta(2+\e))(1-a)/(m_1-m_2)$.
If $F$ satisfies (\ref{RV}) with $1<\alpha<2$, then by Theorem
4.5.2 of Whitt (2000) it is in a normal domain of attraction of
the stable law $S_\alpha(1,1,0)$. Denoting $F_S(dx)$ the
distribution of the stable law $S_\alpha(1,1,0)$ and writing
$a_\alpha(K)=(A/C_\alpha)^{1/\alpha}T^{1/\alpha}(1-\epsilon)^{1/\alpha}$,
we have
\begin{eqnarray*}
\lefteqn{P\left(\frac{\Xi_{[T(1-\epsilon)+1]}-E\sigma
T(1-\epsilon)}{a_\alpha(K)}a_\alpha(K) +\tilde{M}^{(i)}_{[0,\i)}>
(a_i+\kappa)K +m_iT \right)}\\&&= \int_{-\infty}^\infty
P\left(xa_\alpha(K) +\tilde{M}^{(i)}_{[0,\i)}> (a_i+\kappa)K
+m_iT\right)\;dF_S(x)\\&& \geq \int_{[(a_i+\kappa)K
+m_iT]/a_\alpha(K)}^\infty P\left(xa_\alpha(K)
+\tilde{M}^{(i)}_{[0,\i)}> (a_i+\kappa)K +m_iT\right)\;dF_S(x)\\&&
+\int_{+\kappa [(a_i+\kappa)K +m_iT]/a_\alpha(K)}^{\kappa
[(a_i+\kappa)K +m_iT]/a_\alpha(K)} P\left(xa_\alpha(K)
+\tilde{M}^{(i)}_{[0,\i)}> (a_i+\kappa)K +m_iT\right)\;dF_S(x).
\end{eqnarray*}
Recall that
\begin{equation}\label{proprtstab}\overline{F}_S(x)\sim C_\alpha
x^{-\alpha},\quad F_S(-x)={\rm o}(x^{-\alpha}) \quad\mbox{as
$x\to\infty.$}\end{equation} Hence using integration-by-parts or
monotone density theorem, result of Embrechts and Veraverbeke
(1982), for large $K$ we have
\begin{eqnarray*}
\lefteqn{P\left(\frac{\Xi_{[T(1-\epsilon)+1]}-E\sigma
T(1-\epsilon)}{a_\alpha(K)}a_\alpha(K) +\tilde{M}^{(i)}_{[0,\i)}>
(a_i+\kappa)K +m_iT \right)}
\\&&\geq (1-\eta/2)[\overline{F}_S([(a_i+\kappa)K
+m_iT]/a_\alpha(K))\\&&+\frac{C_\alpha
A}{m_i(\alpha-1)}\int_{-\kappa [(a_i+\kappa)K
+m_iT]/a_\alpha(K)}^{\kappa [(a_i+\kappa)K +m_iT]/a_\alpha(K)}
\left((a_i+\kappa)K
+m_iT-xa_\alpha(K)\right)^{-\alpha+1}x^{-\alpha-1}\;dx]
\\&&\geq (1-\eta/2)A[T(1-\e)((a_i+\kappa)K
+m_iT)^{-\alpha}\\&&+
\frac{A}{m_i(\alpha-1)}\left((1+\kappa)\left((a_i+\kappa)K
+m_iT\right)\right)^{-\alpha+1}(1-\overline{F}_S\left(\kappa((a_i+\kappa)K
+m_iT)/a_\alpha(K)\right)\\&&-F_S\left(-\kappa\left((a_i+\kappa)K
+m_iT\right)/a_\alpha(K)\right)) ]\\&&\geq (1-\eta)AT(a_iK
+m_iT)^{-\alpha} + (1-\eta)\frac{A}{m_i(\alpha-1)}\left(a_iK
+m_iT\right)^{-\alpha+1},
\end{eqnarray*}
which completes the proof. \halmos

\begin{Lemma}\label{lowermiesz} Suppose that $F$ satisfies (\ref{RV}).
Given $\eta>0$, it holds for $K$ large enough that
\begin{eqnarray*}
P(M^{(j)}_{N^*} > a_jK, M^{(i)}_{[N^*+1,\i)} > a_iK) \ge (1-\eta)
T \overline{F}(a_jK+m_jT)=(1-\eta)T \overline{F}(a_iK+m_iT).
\end{eqnarray*}
\end{Lemma}
\proof We have
\begin{eqnarray*}
\lefteqn{P(M^{(j)}_{N^*} > a_jK, M^{(i)}_{[N^*+1,\i)} > a_iK)\geq
P(M^{(j)}_{N^*} > a_jK, M^{(i)}_{[N^*+1,\i)} > a_iK,
N^*>T(1-\e))}\\&&\geq P(M^{(j)}_{[T(1-\e)]-1} > a_jK,
M^{(i)}_{[[T(1-\e)],\i)}
> a_iK)-P (N^*<T(1-\e))\\&&\geq
P(M^{(j)}_{[T(1-\e)]-1} > a_jK, S^{(i)}_{[T(1-\e)]} > a_iK) -{\rm
o}(T \overline{F}(a_jK+m_jT))
\\&&\geq
P(M^{(j)}_{[T(1-\e)]-1} > a_jK, \Xi_{[T(1-\e)]} > a_iK +p_iE\zeta
T(1-\e^2), T_{[T(1-\e)]}<E\zeta T(1-\e^2) )\\&&-{\rm o}(T
\overline{F}(a_jK+m_jT))
\\&&\geq
P(M^{(j)}_{[T(1-\e)]-1} > a_jK, \Xi_{[T(1-\e)]}-E\sigma T(1-\e)  >
(1+\e)(a_iK +m_iT)+\e (1-\e)p_iE\zeta
T)\\&&-P(T_{[T(1-\e)]}>E\zeta T(1-\e^2))-{\rm o}(T
\overline{F}(a_jK+m_jT))
\\&&\geq
P(M^{(j)}_{[T(1-\e)]-1} > a_jK, \Xi_{[T(1-\e)]}-E\sigma T(1-\e)  >
(1+\e)(a_jK +m_jT)+\e p_iE\zeta T)\\&&-{\rm o}(T
\overline{F}(a_jK+m_jT))
\\&&\geq
P(M^{(j)}_{[T(1-\e)]-1} > a_jK, S^{(j)}_{[T(1-\e)]}
> (1+\e)a_jK +\e (p_i+p_j)E\zeta T -2\e E\sigma T, T_{[T(1-\e)]}>E\zeta
T(1-\e)^2)\\&& -{\rm o}(T \overline{F}(a_jK+m_jT))
\\&&\geq
P(M^{(j)}_{[T(1-\e)]-1} > a_jK, S^{(j)}_{[T(1-\e)]}
> (1+\e +\xi)a_jK)
-{\rm o}(T \overline{F}(a_jK+m_jT))
\end{eqnarray*}
where $\xi=\e((p_i+p_j)E\zeta -2E\sigma)(1-a)/a_j(m_1-m_2)$.
Denote by $X^{(i)}_k$ the $k$th increment of the random walk
$S_n^{(i)}$. Then \begin{eqnarray*}
\lefteqn{P(M^{(j)}_{[T(1-\e)]-1} > a_jK, S^{(j)}_{[T(1-\e)]}
> (1+\e +\xi)a_jK)}\\&&\geq
P(M^{(j)}_{[T(1-\e)]-1} > a_jK, S^{(j)}_{[T(1-\e)]}
> (1+\e +\xi)a_jK, X_{[T(1-\e)]}<\e a_jK)\\&&\geq
P(M^{(j)}_{[T(1-\e)]-1} > a_jK, S^{(j)}_{[T(1-\e)]-1}
> (1+2\e+\xi)a_jK, X_{[T(1-\e)]}>-\e a_jK)\\&&=
P(S^{(j)}_{[T(1-\e)]-1}
> (1+2\e+\xi)a_jK, X_{[T(1-\e)]}>-\e a_jK)\\&&\geq
P(S^{(j)}_{[T(1-\e)]-1}
> (1+2\e+\xi)a_jK)-P( X_{[T(1-\e)]}<-\e a_jK)\\&&=
P(S^{(j)}_{[T(1-\e)]-1}
> (1+2\e+\xi)a_jK)-{\rm o}(T \overline{F}(a_jK+m_jT)).
\end{eqnarray*}
The assertion of the lemma follows now from the similar
considerations as previously and given in the proof of Lemma
\ref{lowerstad}:
\begin{eqnarray*}\lefteqn{P(S^{(j)}_{[T(1-\e)]-1}
>(1+2\e+\xi)a_jK)}\\&&\geq
P\left((\Xi_{[T(1-\e)]-1}-E\sigma
T(1-\e))/a_\alpha(K)>(1+\e)(1+2\e+\xi)(a_jK+m_jT)/a_\alpha(K)\right)\\&&-P(
T_{[T(1-\e)]-1}>TE\zeta (1+\e)^2)\\&&\geq
\overline{F}_S\left((1+\e)(1+2\e+\xi)(a_jK+m_jT)/a_\alpha(K)\right)
-{\rm o}(T \overline{F}(a_jK+m_jT))\\&&\geq
(1-\eta)\overline{F}_S\left((a_jK+m_jT)/a_\alpha(K)\right)\\&&=
(1-\eta)T \overline{F}(a_jK+m_jT).
\end{eqnarray*}
 \halmos

\subsection{Upper bounds}

\begin{Lemma}\label{upperstad} Suppose that $F$ satisfies (\ref{RV}).
Given $\eta>0$, it holds for $K$ large enough that
\begin{eqnarray*}
P(M^{(i)}_{[N^*+1,\i)} > a_iK) \leq (1+\eta) T
\overline{F}(a_iK+m_iT)+(1+\eta)I^{(i)}_{[T,\i]}(a_iK)
\end{eqnarray*}
\end{Lemma}
\proof Similarly like in the proof of the lower bound given in
Lemma \ref{lowerstad}, for large $K$ and some $\kappa>0$ we have,
\begin{eqnarray*}\lefteqn{P(S^{(i)}_{N^*+1} +\tilde{M}^{(i)}_{[0,\i)}> a_iK)\leq
P(\Xi_{N^*+1} +\tilde{M}^{(i)}_{[0,\i)}> a_iK +p_iT_{N^*+1},
N^*<T(1+\epsilon))}\nonumber\\&&+P(N^*>T(1+\epsilon))\nonumber\\&&\leq
P(\Xi_{[T(1+\epsilon)+1]} +\tilde{M}^{(i)}_{[0,\i)}> a_iK
+p_i\tilde{T})+{\rm o}(I_{[T,\infty]}^{(i)}(a_iK))
\\&&\leq P\left(\frac{\Xi_{[T(1+\epsilon)+1]}-E\sigma
T(1+\epsilon)}{b_\alpha(K)}b_\alpha(K) +\tilde{M}^{(i)}_{[0,\i)}>
(a_i-\kappa)K +m_iT \right)\\&&+{\rm
o}(I_{[T,\infty]}^{(i)}(a_iK))\\
&&\leq (1+\eta/2)[\overline{F}_S([(a_i-\kappa)K
+m_iT]/b_\alpha(K))\\&&+\frac{C_\alpha
A}{m_i(\alpha-1)}\int_{-\kappa [(a_i-\kappa)K
+m_iT]/b_\alpha(K)}^{\kappa [(a_i-\kappa)K +m_iT]/b_\alpha(K)}
\left((a_i-\kappa)K
+m_iT-xb_\alpha(K)\right)^{-\alpha+1}x^{-\alpha-1}\;dx\\&&+\frac{C_\alpha
A}{m_i(\alpha-1)}\int_{\kappa [(a_i-\kappa)K
+m_iT]/b_\alpha(K)}^{[(a_i-\kappa)K +m_iT]/b_\alpha(K)}
\left((a_i-\kappa)K
+m_iT-xb_\alpha(K)\right)^{-\alpha+1}x^{-\alpha-1}\;dx\\&&+
F_S(-\kappa[(a_i-\kappa)K +m_iT]/b_\alpha(K))],
\end{eqnarray*}
where
$b_\alpha(K)=(A/C_\alpha)^{1/\alpha}T^{1/\alpha}(1+\epsilon)^{1/\alpha}$.
The first two increments give required asymptotics since the
second one can be estimated above by
\begin{eqnarray*}\lefteqn{\frac{A}{m_i(\alpha-1)}\left((1-\kappa)\left((a_i-\kappa)K
+m_iT\right)\right)^{-\alpha+1}}\\&&\left(1-\overline{F}_S\left(\kappa((a_i-\kappa)K
+m_iT)/b_\alpha(K)\right) -F_S\left(-\kappa((a_i-\kappa)K
+m_iT)/b_\alpha(K)\right)\right).\end{eqnarray*} The last
increment is ${\rm o}(I^{(i)}_{[T,\i]}(a_iK))$ by
(\ref{proprtstab}). Finally, the third increment could be
estimated in the following way,
\begin{eqnarray*}\lefteqn{\frac{C_\alpha
A}{m_i(\alpha-1)}\int_{\kappa [(a_i-\kappa)K
+m_iT]/b_\alpha(K)}^{[(a_i-\kappa)K +m_iT]/b_\alpha(K)}
\left((a_i-\kappa)K
+m_iT-xb_\alpha(K)\right)^{-\alpha+1}x^{-\alpha-1}\;dx}\\&&=
b_\alpha(K)^\alpha\frac{C_\alpha A}{m_i(\alpha-1)}\int_{\kappa
[(a_i-\kappa)K +m_iT]}^{[(a_i-\kappa)K +m_iT]} \left((a_i-\kappa)K
+m_iT-w\right)^{-\alpha+1}w^{-\alpha-1}\;dw\\&&= ((a_i-\kappa)K
+m_iT)^{-2\alpha+1}T\frac{A^2(1+\e)}{m_i(\alpha-1)}\int_\kappa^1
(1-t)^{-\alpha+1}t^{-\alpha-1}\;dt\\&&={\rm
O}(K^{-2\alpha+2})={\rm o}(K^{-\alpha+1})={\rm
o}(I^{(i)}_{[T,\i]}(a_iK)).
\end{eqnarray*}
 \halmos

\begin{Lemma}\label{uppermiesz} Suppose that $F$ satisfies (\ref{RV}).
Given $\eta>0$, it holds for $K$ large enough that
\begin{eqnarray*}
P(M^{(j)}_{N^*} > a_jK, M^{(i)}_{[N^*+1,\i)} > a_iK) \leq (1+\eta)
T \overline{F}(a_jK+m_jT)=(1+\eta)T \overline{F}(a_iK+m_iT)
\end{eqnarray*}
\end{Lemma}
\proof Similarly like in the proof of the Lemma \ref{lowermiesz}
for large $K$ we have,
\begin{eqnarray*}
\lefteqn{P(M^{(j)}_{N^*} > a_jK, M^{(i)}_{[N^*+1,\i)} >
a_iK)}\\&&\leq  P(M^{(j)}_{[T(1+\e)]-1} > a_jK,
M^{(i)}_{[[T(1+\e)],\i)}
> a_iK)+P (N^*>T(1+\e))
\\&&\leq
P(M^{(j)}_{[T(1+\e)]-1} > a_jK, \Xi_{[T(1+\e)]} > a_iK +p_iE\zeta
T(1-\e^2), T_{[T(1+\e)]}>E\zeta T(1-\e^2) )\\&&+{\rm o}(T
\overline{F}(a_jK+m_jT))
\\&&\leq
P(M^{(j)}_{[T(1+\e)]-1} > a_jK, \Xi_{[T(1+\e)]}-E\sigma T(1+\e)  >
(1-\e)(a_iK +m_iT)-2\e E\sigma T)\\&&+P(T_{[T(1+\e)]}<E\zeta
T(1-\e^2))+{\rm o}(T \overline{F}(a_jK+m_jT))
\\&&\leq
P(M^{(j)}_{[T(1+\e)]-1} > a_jK, S^{(j)}_{[T(1+\e)]}
> (1-\e)a_jK -\e(1+3\e)p_jE\zeta T -2\e E\sigma T,\\&& T_{[T(1+\e)]}<E\zeta
T(1+\e)^2)+P(T_{[T(1+\e)]}>E\zeta T(1+\e)^2)+{\rm o}(T
\overline{F}(a_jK+m_jT))
\\&&\leq
P(M^{(j)}_{[T(1+\e)]-1} > a_jK, S^{(j)}_{[T(1+\e)]}
> (1-\e -\chi)a_jK)
+{\rm o}(T \overline{F}(a_jK+m_jT)),
\end{eqnarray*}
where $\chi=\e[(1+3\e)E\zeta+2E\sigma] (1-a)/a_j(m_1-m_2)$. Then
\begin{eqnarray*}
\lefteqn{P(M^{(j)}_{[T(1+\e)]-1} > a_jK, S^{(j)}_{[T(1+\e)]}
> (1-\e -\chi)a_jK)}\\&&\leq
P(S^{(j)}_{[T(1+\e)]-1}
> (1-2\e-\chi)a_jK)+P(M^{(j)}_{[T(1+\e)]-1} > a_jK, X_{[T(1+\e)]}>\e a_jK)
\\&&\leq
P(S^{(j)}_{[T(1+\e)]-1}
> (1-2\e-\chi)a_jK)+P(M^{(j)}_{[T(1+\e)]-1} > a_jK) P(X_{[T(1+\e)]}>\e a_jK)
\\&&\leq
P(S^{(j)}_{[T(1+\e)]-1}
> (1-2\e-\chi)a_jK)+{\rm O}(K^{-2\alpha +1})
\\&&=
P(S^{(j)}_{[T(1+\e)]-1}
> (1-2\e-\chi)a_jK)+{\rm o}(T \overline{F}(a_jK+m_jT))
\end{eqnarray*}
and
\begin{eqnarray*}
\lefteqn{P(S^{(j)}_{[T(1+\e)]-1}
>(1-2\e-\chi)a_jK)}\\&&\leq
P\left((\Xi_{[T(1+\e)]-1}-E\sigma
T(1+\e))/b_\alpha(K)>(1-\e)(1-2\e\chi)(a_jK+m_jT)/b_\alpha(K)\right)\\&&+P(
T_{[T(1-\e)]-1}<TE\zeta (1-\e^2))\\&&\leq
\overline{F}_S\left((1-\e)(1-2\e-\chi)(a_jK+m_jT)/b_\alpha(K)\right)
+{\rm o}(T \overline{F}(a_jK+m_jT))\\&&\geq
(1-\eta)\overline{F}_S\left((a_jK+m_jT)/b_\alpha(K)\right)\\&&=
(1-\eta)T \overline{F}(a_jK+m_jT),
\end{eqnarray*}
which completes the proof.
 \halmos
\section*{\bf Acknowledgements} { \no We are much indebted to
Florin Avram for his continuing advice on this paper. ZP
acknowledge support by POLONIUM no 09158SD. }

\end{document}